\documentclass[a4paper,11pt]{amsart}
\usepackage[english]{babel}
\usepackage{amsmath,amsfonts,latexsym}
\usepackage{amsthm}
\usepackage{amssymb}
\usepackage[dvips]{graphicx}
\usepackage[dvips]{color}
\usepackage{graphicx}
\usepackage{pgf,pgfarrows,pgfnodes,pgfautomata,pgfheaps}

\usepackage{array}

\usepackage{fancyhdr}
\newtheorem{Thm}{Theorem}[section]
\newtheorem{Lem}[Thm]{Lemma}
\newtheorem{Def}[Thm]{Definition}

\newtheorem{Cla}[Thm]{Claim}
\theoremstyle{remark}

\newtheorem{Rem}[Thm]{\sc Remark}

\numberwithin{equation}{section}

\newcommand{\R}{\mathbb R}
\newcommand{\Z}{\mathbb Z}
\newcommand{\N}{\mathbb N}
\newcommand{\pr}{{\noindent \textsc{Proof}. \, }}  
\newcommand{\epr}{\hspace*{0.1cm}\hfill $\square$\par}

\begin{document}
\title{\bf ON THE VOLUME AND THE NUMBER OF LATTICE POINTS OF SOME SEMIALGEBRAIC SETS}
\author{ {\sf Ha Huy Vui }$^*$}
\address{Institute of mathematics, Vietnam Academy of Sciences, Hanoi, Vietnam \\
18, Hoang Quoc Viet, Cau Giay, Hanoi}
\email{hhvui@math.ac.vn}
\author{ {\sf Tran Gia Loc }$^\dagger$ }
\address{Teacher Training College of Dalat, 29 Yersin, Dalat, Vietnam}
\email{gialoc@gmail.com}
\subjclass[2010]{14B05, 32S20, 34E05, 11H06, 51M20, 52A38, 52A23}
\keywords{Newton polyhedron, the Mikhailov - Gindikin condition, sublevel set, lattice points.}


\maketitle
\vspace{-1cm}
$$\textit{$^*$ Institute of mathematics, Vietnam Academy of Sciences, Hanoi, Vietnam }$$
$$\textit{18, Hoang Quoc Viet, Cau Giay, Hanoi}$$
$$\textit{$^\dagger$ Teacher Training College of Dalat, 29 Yersin road, Dalat, Vietnam}$$

\begin{abstract}
Let $f = (f_1,\ldots,f_m) : \R^n \longrightarrow \R^m$ be a polynomial map; $G^f(r) = \{x\in\R^n : |f_i(x)| \leq r,\ i =1,\ldots, m\}$. We show that
if $f$ satisfies the Mikhailov - Gindikin condition then
\begin{itemize}
\item[(i)] $\text{Volume}\ G^f(r) \asymp r^\theta (\ln r)^k$
\item[(ii)] $\text{Card}\left(G^f(r) \cap \overset{o}{\ \Z^n}\right) \asymp r^{\theta'}(\ln r)^{k'}$, as $r\to \infty$,
\end{itemize}
where the exponents $\theta,\ k,\ \theta',\ k'$ are determined explicitly in terms of the Newton polyhedra of $f$. \\
\indent Moreover, the polynomial maps satisfy the Mikhailov - Gindikin condition form an open subset of the set of polynomial maps having the same Newton polyhedron.
 \end{abstract}
\vskip 0.1 cm
{\bf Keywords and phrases:} Newton polyhedron, the Mikhailov - Gindikin condition, sublevel set, lattice points.
\vskip 0.2 cm

\section{Introduction}
The study of the asymptotic behavior of the volume of sublevel sets and the number of lattice points has attracted a lot of researchers and has found many important applications. In the middle 1970s, A.N. Varchenko and V.A. Vasiliev used Newton polyhedra to study the asymptotic behavior of the volume of sublevel sets and the integrals of real analytic functions in the degenerate situation (\cite{VAR}, \cite{VAS1}, \cite{VAS2}). In particular, sharp estimates for the volume and the integrals were obtained in terms of Newton polyhedra for certain classes of the functions with an isolated minimum at zero.
\vskip 0.1cm
Let $f : \R^n \longrightarrow \R^m$ be a polynomial map. For $r>0$, put 
$$G^f(r) = \{x\in\R^n : |f_i(x)| \leq r,\ i =1,\ldots, m\}, \quad Z^f(r) = G^f(r) \cap \overset{o}{\ \Z^n},$$ 
where $\overset{o}{\ \Z^n} = \{(a_1,\ldots, a_n) \in \Z^n : a_i \neq 0,\ i =1, \ldots, n \}$. Let $\left| G^f(r)\right|$ and $\text{Card}\ Z^f(r)$ be correspondingly the volume of $G^f(r)$ and the cardinal of $Z^{f}(r)$.
\vskip 0.cm
In this paper, we are interested in computing explicitly the exponents arising in the asymptotic formulas for $|G^f(r)|$ and $\text{Card}\ Z^{f}(r)$, as $r \to \infty$.
\vskip 0.1cm
In the case of $m=1$, the asymptotic behavior of $\left| G^f(r)\right|$ plays an important role in many problems of the theory of pseudo-differential operators. \\
\indent The asymptotic behavior of the volume of the set $\{x\in U : |f(x)| < r \}$, as $\ r\to 0$, where $U$ is a small enough neighborhood of a singularity point, concerns the oscillatory integral operators and the scalar oscillatory integrals (see \cite{DNS}, \cite{GRE1}, \cite{GRE2}, \cite{KARO}, \cite{KAR1}, \cite{KAR2}, \cite{P-S}, \cite{PSS2}, \cite{GRE2}, \cite{SEE}).\\
\indent The asymptotic behavior of the volume of the set $\{x\in \R^n : |f(x)| < r \}$, as $\ r\to \infty$, is used in \cite{SIN} to estimate the number of eigenvalues of the Schr$\ddot{o}$dinger operator in $\R^n$.
\vskip 0.1cm
In \cite{DUN}, the asymptotic behavior of $\text{Card}\ Z^f(r)$, where $f$ is a {\it monomial map} was computed and applied in the approximation theory.
\vskip 0.1cm
For the set $G^f(r)$, the following problems arise naturally.
\begin{itemize}
\item[(i)] When are the qualities $\left| G^f(r)\right|$ and $\text{Card}\ Z^{f}(r)$ finite.
\item[(ii)] If they are finite, how to compute the exponents arising in the asymptotic formulas for these qualities?
\end{itemize}
If $f$ is an arbitrary polynomial map then it is very difficult to provide satisfactory answers to these questions, even for the case $n=2$ and $m=1$. However, if the application $f$ satisfies the so called Mikhailov - Gindikin condition, then we can give complete answers to these problems.

\section{Statement of results}
\vskip 0.1cm
For a polynomial $\varphi(x) = \sum a_{\alpha} x^{\alpha} \in \R[x_1,\ldots, x_n]$, we call the support of $\varphi$ the following set 
$$supp(\varphi) := \{ \alpha \in \left(\N \cup \{0\}\right)^n \ : \ a_\alpha \neq 0\}\ .$$
Let $f = (f_1,\ldots,f_m) : \R^n \longrightarrow \R^m$. Put $\Gamma(f) = convex \left( \underset{i=1}{\overset{m}{\cup}} supp(f_i)\right)$, the convex hull of the set $\underset{i=1}{\overset{m}{\cup}} supp(f_i)$. We call $\Gamma(f)$ the {\it Newton polyhedron} of $f$.
\vskip 0.1cm
Let $f_{i}(x) = \sum a^{i}_{\alpha} x^{\alpha}$ and $\Delta$ be a face of $\Gamma(f)$. We put 
$$ f_{i_\Delta}(x) = \sum_{\alpha \in \Delta} a^{i}_{\alpha} x^{\alpha}\ .$$

\begin{Def} We say that $f$ satisfies the Mikhailov-Gindikin condition if for any face $\Delta \subset \Gamma(f)$, we have
$$\max|f_{i_\Delta}(x)| \neq 0, \quad i = 1,\ldots, m ;$$
in $(\R\setminus\{0\})^n$.
\end{Def}

\noindent We denote by $cone \Gamma(f)$ the cone generated by $\Gamma(f)$,  
$$cone \Gamma(f) = \left\{y \ :  \ y = \lambda x \ \text{for} \ \lambda \geq 0 \ \text{and} \  x \in \Gamma(f)\right\}, $$
and by $\Delta^{+}(d)$ the diagonal of the positive orthant in $\R^n$,
$$\Delta^{+}(d) = \{ (d_1,\ldots, d_n) \in \R_{+}^n : d_1 = \ldots = d_n \}.$$ 
Let $D_{\infty}\Gamma(f)$ be the furthest point from the origin in the intersections of the diagonal $\Delta^{+}(d)$ and $\Gamma(f)$, and $\Lambda_{\infty}$ be the face of smallest dimension of $\Gamma(f)$, having $D_{\infty}\Gamma(f)$ as its interior point. 
\vskip 0.1cm
We denote by $k=dim\Lambda_{\infty}$, $D_{\infty}\Gamma(f) = \left(d_{\infty},\ldots, d_{\infty}\right)$, $\theta = \dfrac{1}{d_{\infty}}$ and $v_n = (1,\ldots,1) $.
\vskip 0.1cm
\noindent 
The notation $g \asymp h$ means that there exists positive constants $K_1, K_2$ such that
$$K_1 h \ \leq \  g  \ \leq  \  K_2 h . $$    
 
\begin{Thm}\label{Main1}
Let $f= (f_1,\ldots, f_m): \R^n \longrightarrow \R^m$ be a polynomial map satisfying the Mikhailov-Gindikin condition. Then we have 
\begin{itemize}
\item[(i)] $\left| G^f(r) \right| < \infty$, for any $r > 0$, if and only if   the vector $(1, \ldots, 1)$ belongs to the interior of $cone \Gamma(f)$.
\item[(ii)] If $|G^f(r)| < \infty$, then we have 
$$|G^f(r)| \ \asymp \ r^{\theta} \ln^{n-k-1}r \ , \ \text{as}\quad  r \to \infty \ .$$
\end{itemize}
\end{Thm}

\vskip 0.1cm
Next, we construct the so called {\it complete Newton polyhedron} of $f$. 
\vskip 0.1cm\noindent
For $\alpha, \beta \in \R^n$, we shall write $\alpha \preccurlyeq \beta$, if $\alpha_j \leq \beta_j $ for all $j = 1,\ldots, n$.

\begin{Def}
We call the complete Newton polyhedron of $f$, the polyhedron $\widetilde{\Gamma}(f)$ obtained from $\Gamma(f)$ by adding all the $\alpha \in \overline{\R_{+}^n}$ for which there exists $\beta \in \Gamma(f)$, s.t. $\alpha \preccurlyeq \beta$.
\end{Def}

We denote by $D_{\infty}\widetilde{\Gamma}(f)$ the furthest point from the origin in the intersections of $\Delta^{+}(d)$ and $\widetilde{\Gamma}(f)$. Put  $D_{\infty}\widetilde{\Gamma}(f) = (\widetilde{d}_{\infty}, \ldots, \widetilde{d}_{\infty}) $ and $\theta' =\dfrac{1}{\widetilde{d}_{\infty}}$. Let $\widetilde{\Lambda}_\infty$ be the face having smallest dimension of $\widetilde{\Gamma}(f)$ that contains the point $D_{\infty}\widetilde{\Gamma}(f)$ in its interior. Put $k' = dim\widetilde{\Lambda}_\infty$. 

\begin{Thm}\label{Main2}
Let $f = (f_1,\ldots, f_n) : \R^n \longrightarrow \R^m$ be a polynomial map satisfying the Mikhailov-Gindikin condition. Then we have 
\begin{itemize}
\item[(i)] $Card\ Z^{f}(r) < \infty$ for all positive real numbers $r$, if and only if $cone\Gamma(f) \cap \overset{o}{\ \R^n} \neq \emptyset $.
\item[(ii)] Moreover, if $Card\ Z^{f}(r) < \infty $, we have 
$$Card\ Z^{f}(r)\ \asymp \ r^{\theta'} \ln^{n-k'-1}r  \ \text{as}\  r \to \infty\ .$$
\end{itemize}
\end{Thm}

\begin{Rem}\
\begin{itemize}
\item[(i)] It follows from Theorems 2.2 and 2.4 that, under the Mikhailov - Gindikin condition, the equalities $\theta=\theta',\ k = k'$ hold if and only if $\Lambda = \widetilde{\Lambda}$, i.e $\Lambda$ is a common face of $\Gamma(f)$ and $\widetilde{\Gamma}(f)$.
\item[(ii)] If $f$ is a monomial map, the exponents in the asymptotic formulas for $G^f(r)$ and $Card\ Z^f(r)$ were computed by Dinh Dung in \cite{DUN}. Note that this author has stated his result in terms of  some linear programming problems and did not make use of Newton polyhedra.
\end{itemize}
\end{Rem}
\indent   
Let $\Gamma$ be a convex polytope in $\R^{n}$. Assume that all the vertice of $\Gamma$ belong to $\left(\N \cup \{0\}\right)^n$. 
\vskip 0.1cm
\noindent We define
$$\mathcal{M}_{\Gamma} := \left\{f : \R^n \longrightarrow \R^m : \underset{i=1}{\overset{m}{\cup}} supp(f_i) \subset \Gamma \right\}, \qquad
\mathcal{N}_{\Gamma} := \left\{ f: \R^n\longrightarrow \R^m : \Gamma(f) = \Gamma \right\}\ ,$$
$$\mathcal{D}_{\Gamma} := \left\{ f: \R^n \rightarrow \R^m : \Gamma(f) = \Gamma,\ \text{and}\ f \ \text{satisfies the Mikhailov-Gindikin condition}\right\}\ .$$
\vskip 0.1cm
By the lexicographic ordering in the set of monomials, we can identify $\mathcal{M}_{\Gamma}$ with a finite dimensional space over $\R$, and $\mathcal{N}_{\Gamma}$ and $\mathcal{D}_{\Gamma}$ with subsets in this space.

\begin{Thm}\label{DL_Tapmo}
With the notations above, $\mathcal{D}_{\Gamma}$ is an open subset in $\mathcal{N}_{\Gamma}$, and, consequently, it is an open set in the space $\mathcal{M}_\Gamma$ .
\end{Thm}

\section{Proofs}
Theorem \ref{Main1} and Theorem \ref{Main2} are direct consequences of two following facts 
\begin{itemize}
\item[(i)] Two-side estimation for polynomial functions satisfying the Mikhailov - Gindikin condition.
\item[(ii)] Asymptotic formulas for the volume and the number of lattice points in semi-algebraic sets, defined by monomial inequalities \cite{DUN}. 
\end{itemize}

\vskip 0.1cm 
Let $\varphi(x) = \sum a_{\alpha} x^{\alpha} \in \R[x_1,\ldots, x_n] ,\ \Gamma(\varphi)$ be the Newton polyhedron of $\varphi$ and $V(\varphi)$ be the set of vertices of $\Gamma(\varphi)$.
\vskip 0.1cm
Put $N_{\varphi}(x) = \underset{\alpha \in V(\varphi)}{\sum} |x^{\alpha}|$.

\begin{Thm} (see \cite[p. 204]{GIN}) \label{dinhly_M-G} Two conditions are equivalent
\begin{itemize}
\item[(i)] There is $c > 0$ and $\rho > 0$ such that
$$c N_{\varphi}(x) \ \leq \ |\varphi (x)|,\quad x \in \R^n ,\ |x| > \rho \ .$$
\item[(ii)] For any face $\Delta \subset \Gamma (\varphi)$, and $x \in (\R\setminus \{0\})^n,\ |x| > \rho$, we have
$$\varphi (x) \neq 0, \quad \text{and}\quad \varphi_{\Delta}(x) \neq 0 .$$
\end{itemize}
\end{Thm}

\begin{Rem} It follows from the theorem \ref{dinhly_M-G} and from (\cite[Lemma 1.1]{GIN}) that if $\varphi$ satisfies condition (ii), then exist positive constants $c_1,\ c_2$ and $\rho$ such that 
$$ c_1 N_{\varphi}(x) \ \leq \ |\varphi (x)|  \ \leq  \ c_2 N_{\varphi}(x) \ , \quad x \in \R^n ,\ |x| > \rho \ .$$
\end{Rem}

Now, let us consider the system of monomials 
$$\{ x^{\alpha^1},\ldots, x^{\alpha^s}\},\quad \alpha^i \in \left(\N \cup \{0\}\right)^n  ,\ i=1,\ldots, s\ .$$
For $r>0$, put 
$$G^{\alpha}(r) = \left\{ x \in \R^n : |x|^{\alpha^i} \leq r,\ i=1,\ldots, s \right\}.$$ 
\indent In \cite{DUN}, Dinh Dung computed the first term in asymptotic formulas for volume of $G^{\alpha}(r)$ and for the number of lattice points in $Z^{\alpha}(r) = G^{\alpha}(r)\cap \overset{o}{\ \Z^n}$. We now recall his result.  
\vskip 0.1cm
Consider the following linear programming problem
\begin{equation}\label{BTQHTT_1}
\begin{array}{llll}
x_1 + \ldots + x_n  \ \rightarrow sup \ ;\\
\begin{cases}
\langle x,\alpha^{i}\rangle \ \leq 1 \ , \quad  i = 1,\ldots, s \\
x \in \R^n\ .
\end{cases}
\end{array}
\end{equation}
Let $\theta$ and $M(\alpha)$ be correspondingly the optimal value and the solution set of this problem. Put $p := dim\ M(\alpha)$ and $\Gamma(\alpha) := conv \ \left\{\alpha^1,\ldots, \alpha^s \right\}$, and let $C(\alpha)$ be the cone generated by $\Gamma(\alpha)$. 
\vskip 0.1cm
\begin{Thm} (\cite[Theorem 1]{DUN}) \label{DinhDung1} The volume of $G^{\alpha}(r)$ is finite for all $r>0$ if and only if the vector $(1,\ldots,1)$ belongs to the interior of $C(\alpha)$. Moreover, if volume of  $G^{\alpha}(r)$ is finite for all $r>0$, then
$$\text{volume of} \ G^{\alpha}(r) \ \asymp \ r^{\theta}\ln^{p}r \ .$$
\end{Thm}
\vskip 0.1cm
Next, consider the linear programming problem
\begin{equation}\label{BTQHTT_2}
\begin{array}{llll}
x_1 + \ldots + x_n  \ \rightarrow sup \ ;\\
\begin{cases}
\langle x,\alpha^{i}\rangle \ \leq 1 \ , \quad  i = 1,\ldots, s \\
x \in \R_{+}^n\ .
\end{cases}
\end{array}
\end{equation}
\begin{Thm}(\cite[Theorem 2]{DUN}) \label{DinhDung2} $Card\ Z^{\alpha}(r)$ is finite for any $ r > 0$ if and only if $C(\alpha) \cap \overset{o}{\ \R^n} \neq \emptyset $. Moreover, if this condition is satisfied, then 
$$Card\ Z^{\alpha}(r) \ \asymp \ r^{\theta^{'}} \ln^{p^{'}} r \ ,$$
where $\theta^{'}$ and $p^{'}$ be correspondingly the optimal value and the dimension of the solution set of the  linear programming problem \eqref{BTQHTT_2}.
\end{Thm}
\noindent {\bf Proof of Theorem \ref{Main1}} 
\vskip 0.1cm
Assume that $f : \R^n \longrightarrow \R^m$ satisfies the Mikhailov-Gindikin condition. We put $N_{f}(x) := \underset{\alpha \in V(f)}{\sum} |x|^{\alpha}$, where $V(f)$ is the set of vertices of $\Gamma(f)$. 
\begin{Lem}\label{M-Ksuyrong}
If $f$ satisfies the Mikhailov-Gindikin condition then there exist positive constants $c_1 ,\ c_2 $ and $\rho$ such that
$$c_1 N_{f}(x) \ \leq  \ \max |f_i(x)| \ \leq \ c_2 N_{f}(x)$$
for all $x \in \R^n,\ |x| > \rho$.
\end{Lem}
\pr
Let $ F = \underset{i=1}{\overset{m}{\sum}} f_i^2 $ and $\Gamma(F)$ be the Newton polyhedron of $F$.
It is not difficult to see that if $V(f) = \left\{ \alpha^1, \ldots, \alpha^k \right\}$ is the set of vertices of $\Gamma(f)$ then the set $V(F) = \left\{ 2\alpha^1, \ldots, 2\alpha^k \right\}$ is that of $\Gamma(F)$. In consequence, for every face $\Delta^{'}$ of $\Gamma(F)$, there exists a unique face $\Delta$ of $\Gamma(f)$ such that $\Delta^{'} = 2\Delta$.  
\begin{Cla}\label{Cla_Delta}
Let $\Delta^{'} $ be a face of $\Gamma(F)$ and $\Delta$ be that of $\Gamma(f)$ such that $\Delta^{'} = 2\Delta$. Then 
$$ F_{\Delta^{'}}(x) = \underset{i=1}{\overset{m}{\sum}} f_{i_{\Delta}}^2 (x)\ .$$
\end{Cla}
\pr
\vskip 0.1cm
We begin with a description of a face of a polyhedron in $\R^n$. Let $\Gamma$ be a polyhedron in $\R^n$, $dim \Gamma = n$ and $\Delta$ be its face. Then there exists $q \in \R^n$ such that the restriction of $\langle x, q \rangle$ on $\Gamma$ attains its maximum value at $x$ if and only if $x \in \Delta$. 
\vskip 0.1cm
In fact, if $dim \Delta = n-1$ then $q$ is a normal vector of the hyperplane containing $\Delta$, and $q$ is determined uniquely within a positive factor. If $dim \Delta < n-1$, then $\Delta$ lies on the boundaries of some faces of dimension $n-1$, say $\Delta_1, \ldots, \Delta_l$, where
$$\Delta_i = \left\{ x\in \Gamma : \langle x, q_i \rangle = d(q_i) \right\} $$
and $d(q_i) = \underset{x\in \Gamma}{\sup} \langle x, q \rangle$. Then $\Delta$ can be represented by 
$$\Delta = \left\{ x \in \Gamma : \langle x, q \rangle = d(q) \right\} $$
with $q = \underset{i=1}{\overset{l}{\sum}}t_iq_i  \ , \ \  \underset{i=1}{\overset{l}{\sum}}t_i =1$ and $t_i >0, \ i = 1,\ldots, l$.
\vskip 0.1cm
Now, assume that 
$$\Delta = \left\{ x \in \Gamma(f) : \langle x, q \rangle = d(q) \right\} \quad \text{and} \quad   \Delta^{'} = \left\{ x \in \Gamma(F) : \langle x, q \rangle = 2d(q) \right\}.$$
We write $f_i(x)$ in the form 
$$f_i(x) = h_i(x) + g_i(x)$$
where $h_i(x) = f_{i_{\Delta}}(x)$. In the sum 
$$f_{i}^2(x) = h_i^2(x) + 2h_i(x) g_i(x) + g_i^2(x)$$
every monomial $x^{\alpha}$ satisfying the condition $\langle \alpha, q \rangle = 2 d(q)$, can occur only in the first summand.  Therefore 
$$F_{\Delta^{'}}(x) =  \underset{i=1}{\overset{m}{\sum}} h_i^2 (x) =  \underset{i=1}{\overset{m}{\sum}}f_{i_{\Delta}}^2 (x) .$$
As consequence of this claim, since $f$ satisfies the Mikhailov - Gindikin condition, $F$ satisfies this condition too.
\vskip 0.1cm
By Theorem \ref{dinhly_M-G}, there exist $c >0,\ c' > 0,\ \rho > 0$ such that 
$$|x| > \rho \quad \Rightarrow \quad  c \underset{2\alpha^i \in V(F)}{\sum} \left| x^{2\alpha^i}\right|   \ \leq \  |F(x)| \ \leq \  c'   \underset{2\alpha^i \in V(F)}{\sum} \left| x^{2\alpha^i}\right| $$
for all $x \in \R^n$, where $V(F)$ is the set of vertices of $\Gamma(F)$. And therefore, there exist positive constants $c_1,\ c_2 $ and $\rho_1$ such that 
$$c_1 \underset{\alpha \in V(f)}{\sum} |x^\alpha| \leq max |f_i(x)| \leq c_2 \underset{\alpha \in V(f)}{\sum} |x^\alpha| \ ,$$
for $|x| \geq \rho_1$.
\epr
\vskip 0.1cm
Put 
$$\mathcal{A}(r) := \{x\in\R^n : |x|^\kappa \leq r,\ \kappa\in V(f)\}, \quad \mathcal{B}(r) := \{ x\in \R^n : \max |f_i(x)| \leq r\} .$$
Now, it follows from Lemma \ref{M-Ksuyrong} that there exist constants $\rho_1$ and $\rho_2$ such that 
\begin{equation}\label{bdt_A-B}
\| x\| \geq \rho \quad \Rightarrow \quad \rho_1 |\mathcal{A}(r)| \leq |\mathcal{B}(r)| \leq \rho_2 |\mathcal{A}(r)| .
\end{equation} 
Since 
$$|\mathcal{B}(r)| = \left| \{ |x| \leq \rho : \max|f_i(x)| \leq r \} \cup \{|x| \geq \rho :  \max|f_i(x)| \leq r \}\right|$$
and 
$$\left| \{ |x| \leq \rho :  \max|f_i(x)| \leq r \}\right| \ \leq \ \left| \{ x\in \R^n : \| x\| \leq \rho \}\right|$$
then
$$|\mathcal{B}(r)| \ \asymp \ \left| \{ |x| \geq \rho : \max |f_i(x) | \leq r \}\right|\ , \quad r \to \infty \ .$$
Thus, by \eqref{bdt_A-B}, the proof of Theorem \ref{Main1} is reduced to the problem of computing the exponents in the asymptotic formula of $|\mathcal{A}(r)|$, as $r\to\infty$. For this monomial case, the problem is solved already in \cite{DUN}. \vskip 0.1cm\noindent
Using Theorem \ref{DinhDung1}, we have 
\begin{itemize}
\item[(i)] $|\mathcal{A}(r)| < \infty$ for any $r>0$ if and only if $v_n \in int(cone V(f))$. 
\item[(ii)] If $|\mathcal{A}(r)| < \infty$ then $|\mathcal{A}(r)| \asymp r^{\tilde{\theta}}\ln^{\tilde{k}}r$, where $\tilde{\theta}$ is the optimal value and $\tilde{k}$ is the dimension of the solution set of the following linear programming problem
 \begin{equation}\label{BTQHTT_1'}
\begin{array}{llll}
x_1 + \ldots + x_n  \ \rightarrow sup \ ;\\
\begin{cases}
\langle x,\alpha^{i}\rangle \ \leq 1 \ , \quad  \alpha^i \in V(f) = \{\alpha_1, \ldots, \alpha_s\} \\
x \in \R^n\ ,
\end{cases}
\end{array}
\end{equation}
 \end{itemize}   
To finish the proof of Theorem \ref{Main1}, it rest to prove that $\widetilde{\theta} = \theta$, and $\widetilde{k} = n-k-1$, where the exponents $\theta$ and $k$ are determined in the statement of Theorem \ref{Main1}.
\vskip 0.2cm
We write the linear programming problem above in the form  
\begin{equation}\label{BTQHTT_1''}
\begin{array}{llll}
     \max{\{x_1 + \ldots + x_n\}}, \\
     \begin{cases}
          \alpha^1_{1}x_1 + \ldots + \alpha^1_{n}x_n \leq 1 \\
          \hspace{40pt} \ldots \\
           \alpha^s_{1}x_1 + \ldots + \alpha^s_{n}x_n \leq 1 \\
           (x_1,\ldots,x_n) \in \R^n \ ,
      \end{cases}
\end{array}
\end{equation}
where $\alpha^i  = \left( \alpha^i_1, \ldots , \alpha^i_n \right), \ i=1,\ldots,s $. 
\vskip 0.1cm
Let us consider the dual problem
\begin{equation}\label{BTDN1}
\begin{array}{llll}
     \min{\{u_1 + \ldots + u_s\}}, \\
     \begin{cases}
          \alpha^1_{1}u_1 + \ldots + \alpha^s_{1}u_s = 1 \\
          \hspace{40pt} \ldots \\
           \alpha^1_{n}u_1 + \ldots + \alpha^s_{n}u_s = 1 \\
           u_i \geq 0,\ i = 1,\ldots, s \ .
       \end{cases}
\end{array}
\end{equation}
\noindent The system of linear equations in \eqref{BTDN1} can be rewritten
\begin{equation}\label{PTVT1}
\left(\dfrac{u_1}{\sum_{i=1}^s u_i}\right)\alpha^1 + \ldots + \left(\dfrac{u_s}{\sum_{i=1}^s u_i} \right)\alpha^s = \ \left(\dfrac{1}{\sum_{i=1}^s u_i}, \ldots, \dfrac{1}{\sum_{i=1}^s u_i}\right) \ .
\end{equation}
The point in the left-hand side of \eqref{PTVT1} belongs to $conv\{\alpha^1,\ldots,\alpha^s\}$, and the right-hand side is a point that belongs to $\Delta^+(d)$. On the other hand,  $\dfrac{1}{\sum_{i=1}^s u_i}$ achieves the maximum value when $\underset{i=1}{\overset{s}{\sum}} u_i$ reaches the minimum value. Thus $\underset{i=1}{\overset{s}{\sum}}u_i$  achieves the minimum value at the point $D_{\infty}\Gamma(f) = \left( d_{\infty},\ldots, d_{\infty}\right)$ and $\widetilde{\theta} = \dfrac{1}{d_{\infty}} = \theta$.
\vskip 0.2cm
Put $P := \{ x\in \R^n : \langle x, \alpha^i \rangle \leq 1, \ i =1,\ldots, s\}$ and let $P^{*}$ be the polar set of  $P$, i.e.   
$$ P^{*} = \{ y \in \R^n : \langle x, y \rangle \leq 1,\ \forall x\in P \} .$$
By (\cite{Bro83}, Theorem 9.1, p.57), we have 
$$P^* := conv\{O, \alpha^1,\ldots, \alpha^s \}.$$
According to The Bipolar Theorem (\cite{Bar02}), $(P^*)^* = P$.
\vskip 0.2cm
Let $\Lambda_{max}$ be the solution set of the problem \eqref{BTQHTT_1'}. Then $\Lambda_{max}$ is a face of $P$. Put 
 $$ \Lambda_\infty = \{y \in P^* : \langle x,y \rangle =1,\ \forall x\in \Lambda_{max} \}.$$
Then, $\Lambda_\infty$ is a face of $P^{*}$ and $\Lambda_{\infty}^{*}$, the polar set of $\Lambda_\infty$, is equal to $\Lambda_{max}$. We see that if $x\in \Lambda_{max}$, then $x_1 + \ldots +x_n = \theta$. Therefore $\left\langle D_\infty\Gamma(f), x \right\rangle = \dfrac{1}{\theta}\left( x_1 + \ldots + x_n \right) = 1$, hence $D_{\infty}\Gamma(f) \in \Lambda_\infty$. Since $\Lambda_\infty$ does not contain the origin, $\Lambda_\infty$ is the face of $\Gamma(f)$ containing the point $D_{\infty}\Gamma(f)$.
\vskip 0.2cm
Now, since $\Lambda_\infty = \Lambda^{*}_{max}$, we have 
$$dim\ \Lambda_{max} = \tilde{k} =n- dim\ \Lambda_\infty - 1 = n-k-1 .$$  
\epr 
 
\subsection{Proof of Theorem \ref{Main2}}\
\vskip 0.1cm
As in the proof of Theorem \ref{Main1}, the proof of  Theorem \ref{Main2} is reduced to the problem of computing the exponents in the asymptotic formula $card\ Z^{V(f)}(r)$, as $r\to\infty$, where $card\ Z^{V(f)}(r) = \mathcal{A}(r) \cap \overset{o}{\ \Z^n}$. Using Theorem \ref{DinhDung2}, we have
\begin{itemize}
\item[(i)] $Card\ Z^{V(f)}(r)  < \infty$ for any $r>0$ if and only if $cone V(f) \cap \overset{o}{\ \R^n} \neq \emptyset $.  
\item[(ii)] If $Card\ Z^{V(f)}(r)  < \infty$ then $Card\ Z^{V(f)}(r) \asymp r^{\tilde{\theta}'}\ln^{\tilde{k}'}r$, where $\tilde{\theta}'$ is the optimal value and $\tilde{k}'$ is the dimension of the solution set of the following linear programming problem
 \begin{equation}\label{BTQHTT_2'}
\begin{array}{llll}
x_1 + \ldots + x_n  \ \rightarrow sup \ ;\\
\begin{cases}
\langle x,\alpha^{i}\rangle \ \leq 1 \ , \quad  \alpha^i \in V(f) = \{\alpha_1, \ldots, \alpha_s\} \\
x \in \R_{+}^n\ ,
\end{cases}
\end{array}
\end{equation}
 \end{itemize}   
We will show that $\widetilde{\theta}' = \theta'$, and $\widetilde{k}' = n -k' -1$, where the exponents $\theta'$ and $k'$ are determined in the statement of Theorem \ref{Main2}.   
\vskip 0.2cm
We write the linear programming problem \ref{BTQHTT_2'} in the form  
\begin{equation}\label{BTQHTT_2''}
\begin{array}{llll}
     \max{\{x_1 + \ldots + x_n\}}, \\
     \begin{cases}
          \alpha^1_{1}x_1 + \ldots + \alpha^1_{n}x_n \leq 1 \\
          \hspace{40pt} \ldots \\
           \alpha^s_{1}x_1 + \ldots + \alpha^s_{n}x_n \leq 1 \\
           (x_1,\ldots,x_n) \in \R^n,\ x_j \geq 0,\ j=1,\ldots, n ,
      \end{cases}
\end{array}
\end{equation}
where $\alpha^i =(\alpha^i_1,\ldots,\alpha^i_n),\ i = 1,\ldots,s$.
\vskip 0.1cm\noindent 
Let us consider the dual problem  
\begin{equation}\label{BTDN2}
\begin{array}{llll}
     \min{\{u_1 + \ldots + u_s\}}, \\
     \begin{cases}
          \alpha^1_{1}u_1 + \ldots + \alpha^s_{1}u_s \geq 1 \\
          \hspace{40pt} \ldots \\
           \alpha^1_{n}u_1 + \ldots + \alpha^s_{n}u_s \geq 1 \\
           u_i \geq 0,\ i = 1,\ldots, s \ .
       \end{cases}
\end{array}
\end{equation}
The system of linear inequations in \eqref{BTDN2} can be rewritten 
\begin{equation}\label{BPTVT1}
\left(\dfrac{u_1}{\sum_{i=1}^s u_i}\right)\alpha^1 + \ldots + \left(\dfrac{u_s}{\sum_{i=1}^s u_i} \right)\alpha^s \ \geq \ \left(\dfrac{1}{\sum_{i=1}^s u_i}, \ldots, \dfrac{1}{\sum_{i=1}^s u_i}\right) \ .
\end{equation}
Put $\gamma^l = \left(\dfrac{u_1}{\sum_{i=1}^s u_i}\right)\alpha^1 + \ldots + \left(\dfrac{u_s}{\sum_{i=1}^s u_i} \right)\alpha^s$ and $\gamma^r = \left(\dfrac{1}{\sum_{i=1}^s u_i}, \ldots, \dfrac{1}{\sum_{i=1}^s u_i}\right)$. Since $\gamma^l \in \Gamma(f)$ and $\gamma^r \leq \gamma^l$, we have $\gamma^r \in \widetilde{\Gamma}(f)$, the complete Newton polyhedron of $f$. 
\\
On the other hand, $\dfrac{1}{\sum_{i=1}^s u_i}$ achieves the maximum value when $\underset{i=1}{\overset{s}{\sum}}u_i$ reaches the minimum value. It follows that $\underset{i=1}{\overset{s}{\sum}}u_i$ reaches the minimum value at the point $D_{\infty}\widetilde{\Gamma}(f) = \left(\widetilde{d}_\infty, \ldots, \widetilde{d}_\infty \right)$ and $\widetilde{\theta}' = \dfrac{1}{\widetilde{d}_\infty}= \theta'$. 
\vskip 0.2cm
Put $\widetilde{P} = \{ x \in \overline{\R^n_+} : \langle x, \alpha^i \rangle \leq 1, \ i =1,\ldots, s\}$. Then, $\widetilde{P}$ is a bounded convex polyhedron having faces which intersect the axes $Ox_j$ at the points $A_j,\ j=1,\ldots, n$, and containing the origin $O$. Hence, 
$$\widetilde{P}= conv\{O, A_1,\ldots, A_n, \alpha_1,\ldots, \alpha_s\} .$$
From the properties of polar sets (see \cite{Gr03, Zie95}) we have, $\widetilde{P}^* = \underset{\beta\in V(\widetilde{P})}{\bigcap} K(\beta,1)$, where $K(\beta,1) = \{x\in\R^n :  \langle \beta, x \rangle \leq 1 \}$ and $V(\widetilde{P})$ is the set of the vertices of $\widetilde{P}$.
Hence, 
$$\widetilde{P}^* \cap \overline{\R^n_+} = \widetilde{\Gamma}(f) .$$
Let $\widetilde{\Lambda}_{max}$ be the solution set of the problems \eqref{BTQHTT_2'}. Put 
$$\widetilde{\Lambda}_\infty := \{ y\in \widetilde{P}^* : \langle x,y \rangle = 1,\ \text{for all }\  x\in \widetilde{P}\}. $$ 
Then, $\widetilde{\Lambda}_\infty$ is a face of $\widetilde{P}^*$ and $(\widetilde{\Lambda}_{max})^* = \widetilde{\Lambda}_\infty$. We see that if  $x\in\widetilde{\Lambda}_{max}$, then $x_1 + \ldots + x_n = \theta'$. Therefore 
$\left\langle D_\infty\widetilde{\Gamma}(f), x \right\rangle = \dfrac{1}{\theta'}(x_1 + \ldots + x_n) = 1$. Hence, $D_\infty \widetilde{\Gamma}(f) \in \widetilde{\Lambda}_\infty$. Since $\widetilde{\Lambda}_\infty$ does not contain the origin, $\widetilde{\Lambda}_\infty$ is the face of $\widetilde{\Gamma}(f)$, which contains the point $D_\infty\widetilde{\Gamma}(f)$ and      
$$dim \widetilde{\Lambda}_{max}  = \tilde{k'} = n - dim \widetilde{\Lambda}_\infty - 1 = n - k' - 1 .$$ 

\subsection{Proof of Theorem \ref{DL_Tapmo}}\
\vskip 0.1cm
Put $\Omega := 2\Gamma$ and $N_{\Omega} := \left\{ h \in \R[x_1,\ldots,x_n] : \Gamma(h) = \Omega \right\}$. We consider the map
$$ F : \mathcal{N}_{\Gamma} \longrightarrow N_{\Omega},\quad g=(g_1,\ldots,g_m) \longmapsto F_g = \sum_{i=1}^m g_i^2\ ,$$
where $F_g(x) = \underset{i=1}{\overset{m}{\sum}} g_i^2(x)$. It is obvious that $F$ is a continuous mapping. Put
$$A_{\Omega} := \Big\{ f\in N_{\Omega} : \text{there exist}\  r > 0,\ c > 0\ \text{such that}\ \|x\| \geq r \Rightarrow f(x) \geq c\underset{\alpha\in V(\Omega)}{\sum} x^{\alpha}\Big\}. $$
\begin{Cla}\label{Cla1}
$g=(g_1,\ldots,g_m)\in \mathcal{D}_{\Gamma}$ if and only if $F_g \in A_{\Omega}$.
\end{Cla}
\pr
Let $g=(g_1,\ldots,g_m) \in \mathcal{D}_{\Gamma}$. Then $\Gamma(g) = \Gamma$, and for any face $\Delta$ of $\Gamma(g)$, we have
$$\max|g_{i,\Delta}(x)| \neq 0 \ \text{for all}\  x\in (\R^n\setminus \{0\})^n.$$
Let $\Delta' $ be a face of $\Omega$, $\Delta ' = 2\Delta$. By Claim \ref{Cla_Delta}, we have
$$ F_g(x) = \sum_{i=1}^m g_i^2(x) \neq 0,\ \text{and}\ \left(F_g\right)_{\Delta'}(x) = \sum_{i=1}^m g_{i,\Delta}^2(x) \neq 0.$$
Therefore $F_g$ satisfies the Mikhailov - Gindikin. By Theorem \ref{dinhly_M-G}, there exists $c > 0$ and $\rho > 0$ such that
$$| F_g(x)| \geq c \sum_{\alpha\in V_{\Omega}} |x^{\alpha}| = c \sum_{\alpha\in V_{\Omega}} x^{\alpha},\ \text{for all} \ x\in\R^n \ \text{satisfying} \ |x| > \rho \ .$$
since all the point $\alpha \in V_{\Omega}$ have even coordinates.
\vskip 0.1cm
We now prove the converse.\\
Let $F_g \in A_{\Omega}$. Then $\Gamma(F_g) = 2\Gamma(g)$ and there exist the numbers $r>0,\ c>0$ such that
$$\|x\| \geq r \Longrightarrow F_g(x) = \sum_{i=1}^m g_i^2(x) \geq c\sum_{\alpha\in V_{\Omega}} x^{\alpha} .$$
Let $\Delta$ be a face of $\Gamma(g)$, and $\Delta ' = 2\Delta$ be the corresponding face of $\Omega$. Let $q = (\rho_1,\ldots,\rho_n)$ be an interior point of the normal cone of $\Delta$. Then
$$\Delta = \{x\in \Gamma : \langle x,q \rangle = d(q) \}, $$
and $\langle x,q \rangle < d(q)$ if $x\in \Gamma \setminus \Delta$.
\vskip 0.1cm
\noindent Take a point $x^0 \in (\R\setminus \{0\})^n$, we see that
$$
\begin{array}{lll}
F_g(t^{\rho_1}x_1^0,\ldots, t^{\rho_n}x_n^0) &=& \left(F_g\right)_{\Delta'} (t^{\rho_1}x_1^0,\ldots, t^{\rho_n}x_n^0) + \text{lower order terms in t}\\
\\
                                                                        & = & t^{2d(q)}\left( F_g \right)_{\Delta'}(x^0) + o\left(t^{2d(q)}\right),\quad t \to \infty.
\end{array}
$$
Hence
$$ t^{2d(q)}\left( F_g \right)_{\Delta'}(x^0) + o\left(t^{2d(q)}\right) \geq c\sum_{\alpha\in V_{\Omega}} x^{\alpha},\quad \text{for t sufficiently large}.$$
On the other hand, each point of $V_{\Omega}$ has even coordinates, this inequality implies that
$$\left( F_g \right)_{\Delta'} (x^0) > 0.$$
By Claim \ref{Cla_Delta},
$$\left( F_g \right)_{\Delta'} (x^0) = \sum_{i=1}^m g_{i,\Delta}^2 (x^0) > 0.$$
It follows that,
$$\max|g_{i,\Delta}(x^0)| \neq 0.$$
\epr

\begin{Cla}\label{Cla2}
$A_{\Omega}$ is an open subset of $N_{\Omega}$.
\end{Cla}

\pr
Assume $f(x) = \sum c_{\alpha}x^{\alpha} \in A_{\Omega}$. Therefore, there exist $r>0,\ c>0$ such that
$$\|x\| \geq r \ \Longrightarrow \ f(x) \geq c \sum_{\alpha\in V_{\Omega}} x^{\alpha}\ .$$
We shall show that there exists a number $\epsilon >0$ such that, if $|\delta_{\alpha}| < \epsilon$, for any $\alpha \in \Omega$, then $\widetilde{f}(x) := \underset{\alpha\in\Omega}{\sum} (c_{\alpha} + \delta_{\alpha})x^{\alpha} \in A_{\Omega}$.
\vskip 0.1cm
In fact, if $\|x\| \geq r$ then
\begin{equation}\label{BDT6}
\widetilde{f}(x) \geq \sum_{\alpha\in\Omega} c_{\alpha} x^{\alpha} - \sum_{\alpha\in\Omega} |\delta_{\alpha}||x|^{\alpha}\ .
\end{equation}
By \cite{G-V} (Lemma 1, p. 160), if $\alpha \in \Omega \cap \N^n$ then
\begin{equation}\label{BDT7}
|x|^{\alpha} \leq \sum_{\alpha\in V_{\Omega}}x^{\alpha}.
\end{equation}
Thus
\begin{equation}\label{BDT8}
\sum_{\alpha \in \Omega \cap \N^n} |\delta_{\alpha}||x|^{\alpha} \leq    \sum_{\alpha \in \Omega \cap \N^n}\epsilon |x|^{\alpha} \leq \epsilon \eta \sum_{\alpha\in V_{\Omega}}x^{\alpha},
\end{equation}
where $\eta$ is the number of integer points in $\Omega$.
\vskip 0.1cm
Combining the inequalities \eqref{BDT6}, \eqref{BDT7}, and \eqref{BDT8} we get the following inequality
$$\widetilde{f} (x) \geq  (c-\epsilon\eta) \sum_{\alpha\in V_{\Omega}} x^{\alpha}, $$
for all $x\in \R^n$ satisfying $\|x\| \geq r$.
\vskip 0.1cm
Thus, if $\epsilon = \dfrac{c}{2\eta}$ then $\widetilde{f}(x) \geq \dfrac{c}{2} \underset{\alpha\in V_{\Omega}}{\sum} x^{\alpha}$. Therefore $\widetilde{f}(x) \in A_{\Omega}$ and the claim \ref{Cla2} is proved.
\epr
\noindent  We continue the proof of Theorem \ref{DL_Tapmo}.
\vskip 0.1cm
Assume that $g_{0}\in \mathcal{D}_{\Gamma}$. We will show that there exists an open neighborhood $U(g_0)$, s.t. $U(g_0) \subset   \mathcal{D}_{\Gamma}$. By the claim \ref{Cla1}, since $g_{0}\in \mathcal{D}_{\Gamma}$ we have 
$$ F(g_{0}) \in A_{\Omega}. $$
By Claim \ref{Cla2}, there exist an open set $V \subset A_{\Omega}$, containing $F(g_0)$. The mapping $F : \mathcal{N}_{\Gamma} \longrightarrow N_{\Omega}$ is continuous, then there exists an open neighborhood $U(g_0)$ of $g_0$, such that 
$$F\left( U(g_0) \right) \subset V .$$
Takes any element $g\in U(g_0)$, we have $F(g)\in V$. Hence $F(g) \in A_{\Omega}$. Now, it follows from Claim \ref{Cla1}, $g \in \mathcal{D}_{\Gamma}$. Thus the open set $U(g_0)$ is contained $\mathcal{D}_{\Gamma}$. 
\epr

\vskip 1cm
\indent {\bf Acknowledgments}. 
\vskip 0.1cm
We would like to thank Professor L\^e Dung Trang, who provided insight and expertise that greatly assisted this research.
\vskip 0.1cm
This paper is supported by Vietnam’s National Foundation for Science and Technology Development (NAFOSTED).

\end{document}